\documentclass[11pt]{amsart}
\usepackage{amssymb,amsmath,latexsym,amsthm}

\usepackage{times}
\usepackage{dsfont}
\usepackage[all]{xy}

 \usepackage[english,francais]{babel}

 \usepackage{url} 
 \usepackage{hyperref}

\input amssym.def
\newtheorem{theorem}{{Theorem}}[section]
\newtheorem{isom.ext}[theorem]{{Trivial isometric extension}}%[section]
\newtheorem{definition}[theorem]{{Definition}}%[section]
\newtheorem{remark}[theorem]{{Remark}}%[section]
\newtheorem{remarks}[theorem]{{Remarks}}%[section]
\def\C{\mathbb{C}}
\def\R{\mathbb{R}}

\def\Z{\mathbb{Z}}

\def\T{\mathbb{T}}

\def\GL{{\sf{GL}}\,}

\def\SL{{\sf{SL}}\,}

\def\Is{{\sf{Iso}}}
\def\Diff{{\sf{Diff}}}

\def\det{{\sf{det}}}

\begin{document}

\selectlanguage{english}
%***

%***

\title{Singular Riemannian metrics, sub-rigidity vs rigidity} 
%Rigidity of degenerate metrics, and further remarks}
\author[S. Bekkara]{Samir Bekkara$^\star$}
\address{S. Bekkara, UST-Oran, Algeria  }
\email{samir.bekkara@gmail.com}

\thanks{${}^\star$ Supported in part by  the ANR Geodycos of the ENS-Lyon and the project RIAMI of the CIMPA}
 
\author[A. Zeghib]{Abdelghani Zeghib}
\address{A. Zeghib, 
CNRS, UMPA, ENS-Lyon, France  
}
\email{zeghib@umpa.ens-lyon.fr \hfill\break\indent
 \url{http://www.umpa.ens-lyon.fr/~zeghib/}}
\date{\today}

\begin{abstract} We analyze sub-Riemannian  and lightlike metrics 
%the case of two cases of
from the point of view of their rigidity as  geometric structures.  Following Cartan's and Gromov's  formal definitions, they are never rigid,  yet, in generic cases, they naturally give rise   to rigid geometric structures!?

\end{abstract}

\maketitle

%\section*

\tableofcontents
 
\section{On sub-Riemannian metrics}

The following are variations on the concept of rigidity of geometric structures in 
a somehow ``paradoxical'' situation:

%\footnote{Detailed proofs of the statements in this note will appear in \cite{BZ}}

%We will deal here with variations around the following  somehow ``paradoxical'' situation%\footnote{Detailed proofs of the statements in this note will appear in \cite{BZ}}:

 \subsubsection*{\sc Sub-Riemannian metrics}
   A 
 {\bf sub-Riemannian} structure $(M, D, h)$ consists in giving on a manifold $M$ a hyperplane field $D \subset TM$
 together with a metric $h$ defined on $D$ (and thought of as infinite on $TM-D$). An isometry
 of $h$ is a diffeomorphism preserving the structure.   
 
%Consider a contact sub-Riemannian structure: that is a Riemannian metric $h$ defined on  a %contact hyperplane field $D $ on a manifold $M^{2d+1}$. 

The   hyperplane field may be defined locally as the kernel of a 1-form $\omega_0$. There is however no canonical choice, any form $\omega = f \omega_0$ defines the same hyperplane.

\subsubsection*{Integrable case}If $D$ is integrable, say it defines a foliation ${\mathcal F}$, then $h$ is nothing but a leafwise  Riemannian metric. We have for instance   the particular global product case: $M = N \times S$,  where the leaves $N \times \{*\}$ are endowed with a same metric $h_0$. Any family $(f_s)_{s\in S}$ in $\Is(N, h_0)$ determines an isometry of 
$(M,D, h)$.

%$h$ consists of a family of Riemannian metrics $\{h_s\}_{s \in S}$ on $N$. Any map $f: (n, s)  %\to f_s(n), s)$, where $f_s \in \Is (N, h_s)$ is an isometry of $h$. In particular, if all %the $h_s$ are the same, then any 

%\subsubsection*{Contact case} 

\subsubsection{A correspondence: contact sub-Riemannian $\to$ Riemannian} \label{h.bar}

We will henceforth assume  that $D$ is a contact hyperplane, i.e.  $\omega_0 \wedge (d\omega_0)^d $ is a volume form where $\dim M = n = 2 d +1$. Let us recall the (classical) construction of a natural Riemannian metric $\bar{h}$ associated to $(D, h)$. 
Observe that   $d\omega_0$ is a symplectic form on $D$ and that for a function $f$,   $d (f\omega_0) = f d\omega_0$ (on $D$). Assume $D$  orientable, and let $\alpha$ be the Riemannian volume form derived from $h$ on it. Writing that $d\omega^d = \alpha$ on $D$ determines uniquely $f$, in other words the  Riemannian metric (together with an orientation) allows one  to choose a canonical  contact form, say $\omega_0$. Let $R$ be the Reeb field of $\omega_0$: $i_R d\omega_0 = 0$, and $\omega_0(R) = 1$. We extend the metric $h$ to a metric $\bar{h}$ on $TM$ by declaring  that $R$ is  unit and orthogonal to $D$. The orientation is actually irrelevant since the inverse orientation gives the same metric. 

\subsubsection*{Isometry groups of Lie type} Summarizing up,  a contact sub-Riemannian metric generates a Riemannian metric. In particular the isometry group  $\Is(M, h)$ is a (closed) sub-group of $\Is(M, \bar{h})$. Similarly 
for pseudo-groups of local isometries (i.e. isometries defined between open sets of $M$, and composed when this is possible). It then follows that the 
isometry group of $(M, h)$ as well as its local isometry pseudo-group  are of Lie type (of finite dimension). 

\subsubsection*{\sc Cartan's finite type condition} (see \cite{Kob}) Let $H$ be a subgroup of $\GL_n(\R)$. An $H$ structure on a 
manifold $M^n$ is a reduction of the structural group of its frame bundle $\GL^{(1)} (M)= P(M)$ to $H$. Equivalently, this is a section of $P(M)/H$ (assume here to simplify that $H$ is closed). A Riemannian metric corresponds to a $O(n)$-structure. A sub-Riemannian metric corresponds to an $H$ structure where $H$ is the subgroup of $\GL_n(\R)$ preserving $\R^{n-1}$ and the standard Euclidean product on it. Its elements have the form:
$$\left(
\begin{array}{cc}
A & \overrightarrow{u} \\
0 & b
\end{array}
\right)$$
where $ A \in O(n-1),  \overrightarrow{u} \in \R^{n-1}$ and $b \in \R$.

Following Cartan, one associates to $H$ its $k$-prolongation $\mathcal H_k$,  a space of symmetric $(k+1)$-multi-linear forms 
on $\R^n$ with values in $\R^n$. If $A \in \mathcal H_k$, and $v_1, \ldots, v_k \in \R^n$ 
are fixed, then $v \to A(v, v_1, \ldots, v_k)$ belongs to $End(\R^n)$. By definition of a prolongation, the last element is assumed to be in the Lie algebra of $H$.

An $H$-structure has a {\bf finite type} $k \in \Bbb N$,  if $\mathcal H_k = 0$. The principal result  of Cartan theory is that the isometry group of an $H$-structure of finite type is a Lie group. The remarkable fact here is that being of finite type depends only on $H$ (as a subgroup of $\GL_n(\R)$)  and not on the structure itself. As an example, a sub-Riemannian metric has infinite type, no matter it is 
integrable or contact! If fact, the test of finiteness of type  concerns the case of the flat translation-invariant $H$-structure on $\R^n$. The flow of a vector field $V$ preserves this structure iff, for any $x \in \R^n$,  the derivative $D_xV $ belongs to the Lie algebra $\mathcal H$.  The $(k+1)$-coefficient in the Taylor development in a linear coordinates of $V$ belongs to $\mathcal H_k$. Hence, finite type means $V$ is polynomial.

\subsubsection*{\sc Gromov's rigidity}  (see \cite{Bal, Be, Can1, D-G, Gro}) Gromov's definition of geometric structures consists essentially in giving up the ``infinitesimal homogeneity'' in the case of Cartan's $H$-structures. As examples, functions, vectors fields are geometric structures, and also is  a `` finite union'' of geometric structures. Isometries are defined naturally.  Gromov introduces a  rigidity condition which coincides with finiteness of type in the case of $H$-structures. 

\subsubsection*{Rough-definition} If $\sigma$ is such a structure on a manifold $M$, and $x \in M$, let 
$\Is^{Loc}_x(\sigma)$ be the group of (germs) of isometries defined in a neighborhood of $x$ and fixing $x$.  For an integer $k$,   denote by  $\Diff^k_x(M)$ the group of $k$-jets of diffeomorphisms of $M$ fixing $x$. We have a map $jet^k_x: \Is^{Loc}_x(\sigma) \to \Diff^k_x(M)$. The intuitive idea of rigidity (of order $k$) is that $jet^k_x$ is injective: an isometry is fully determined by its $k$-jet. We say in this case that $\sigma$ is {\bf Iso-rigid at order $k$}. For example, a Riemannian metric is Iso-rigid at order $1$: an isometry is determined by giving its derivative at some point. In the case of sub-Riemannian metrics, $jet^1_x$ is injective in the contact case (since it generates a Riemannian metric), but for no $k$,  $jet^k_x$ is   injective in the integrable one. 

\subsubsection*{Definition } We then conclude a divergence between this intuitive formalization of rigidity and Cartan's finiteness of type. The true Gromov's definition is actually of infinitesimal nature. For an $H$-structure $\sigma$, one defines the group $\Is^{k+1}_x(\sigma) \subset \Diff^{k+1}_x(M)$ as the group 
of $(k+1)$-jets of diffeomorphisms preserving $\sigma$ up to order $k$ at $x$. 
 For example, if $\sigma$ is a Riemannian metric, then a (local)-diffeomorphism $f$ gives rise to a $(k+1)$-infinitesimal isometry $\in \Is^{k+1}_x(\sigma)$ if $f^*\sigma- \sigma$ vanishes up to order $k$ at $x$. 
(In the general case of  a geometric structure $\sigma$ of order $r$,   $f$ is a $(k+r)$-isometry if $f^*\sigma$ and $\sigma$ have the same $k$-jet at $x$).
 The true definition of {\bf $k$-rigidity} is that $jet^k: \Is^{k+1}_x(\sigma)  \to \Is^k_x(\sigma)$ is injective for any $x$. 

\subsubsection*{Example} Let us see how this injectivity default happens in the example   of the  contact form $\omega = dz +xdy -ydx$ on $\R^3$, endowed with the restriction of $dx^2 + dy^2$. It corresponds to a left invariant contact sub-Riemannian structure on the Heisenberg group, and hence it  is homogeneous. Consider $f: (x, y, z) \to (x +\delta (z), y + \delta (z), z)$. Assume $\delta(0) = 0$, then, $f(0) = 0$.  Thus, $jet^{k+1}_0f$ determines a $(k+1)$-isometry at 0, iff, $\frac{ \partial \delta}{\partial z^k}(0) = 0$ or equivalently 
$jet^k_0(f) = 1$  (that is $f$ has the same $jet$ as the identity). So any such $\delta$ with a non-trivial $\frac{ \partial \delta}{\partial z^{k+1}}(0) \neq  0$ determines an isometry violating the injectivity of $\Is^{k+1} \to \Is^k$. 

\subsubsection*{Remark}  Gromov's definition strictly coincides with  finiteness of type in the case of $H$-structure (see \cite{Bal}, Example 3.17), and thus sub-Riemannian metrics are $k$-rigid for no $k$.

Let us end this criticism   on the rough definition Iso-rigidity     by noting that  a generic geometric structure (e.g. a Riemannian metric) has no non-trivial local isometries, in which case the local rigidity condition is empty. In contrast, it is the infinitesimal rigidity condition (even empty) that allows one to associate ``rigid'' (solid!) objects to $\sigma$, independently of the fact that it has or not local isometries. Indeed, it is proved in both Cartan and Gromov situations, that $k$-rigidity (or $k$-finiteness of type) allows one to construct a parallelism canonically  associated to $\sigma$ 
defined on the $k$-frame bundle $\GL^{(k)}(M) \to M$ (this is the usual frame bundle for $k= 1$) (see for instance \cite{Bal} and Theorem 2 in \cite{Can1}). 
This produces in particular the Levi-Civita connection and hence geodesics for Riemannian metrics.   Also, the Lie group property is proved by  means of this framing.

% Let $\Phi$ be a geometric structure on  a differentiable manifold $M$.

% One is intersetd The common feeling is that 
 
 %We all have

\section{Lightlike   metrics}

\subsubsection*{Duality} Our original motivation was to study rigidity of lightlike metrics. They are simply $H$-structures where $H$ consists of matrices:    
%$O(0, 1, n)$
$$\left(
\begin{array}{cc}
b & \overrightarrow{u} \\
0 & A
\end{array}
\right)$$
where $ A \in O(n-1),  \overrightarrow{u} \in \R^{n-1}$ and $b \in \R$.
Observe that this  is exactly the dual of the group  defining sub-Riemannian metrics,  that is, 
the automorphism $A \in \GL_n(\R) \to {A^*}^{-1} \in \GL_n(\R)$ sends one group onto  the other.

More geometrically, one defines a {\bf lightlike} scalar product  on  a vector space   as a positive 
symmetric bilinear form having a kernel of dimension one. A lightlike metric on a manifold $M$ is a tensor which is a lighlike scalar product on each tangent space. More generaly, a lighlike metric on a vector bundle $E \to M$ consists in giving a 1-dimensionnal sub-bundle 
$N \subset TM$ together with a Riemannian metric on $E/N$. If one defines a sub-Riemannian metric on $E\to M$ as a codimension 1 sub-bundle $D \subset E$ endowed with a Riemannian metric, then one gets a duality: 
\begin{center}
lightlike metric on $E $ $\leftrightarrow$ sub-Riemannian metric on (the dual) $E^*$.
\end{center}
In other words, a lightlike metric $g$ on a manifold 
$M$ consists in giving a line sub-bundle (direction field) $ N \subset TM$, and a Riemannian metric 
on $TM/N$. The direction field $N$ and the 1-dimensional foliation $\mathcal N$ that it generates are called {\bf characteristic}.

\subsection*{Natural situations} Lightlike metrics appear naturally as induced metrics on submanifolds of Lorentz manifolds. Indeed let $(L, h)$ be a Lorentz manifold, and $M \subset L$ a submanifold
such that for any $x \in M$, the restriction $h_x$ on $T_xM$ is degenerate.
%(i.e. not a pseudo-Riemannian metric). 
Then, this is a lightlike metric on $M$, i.e. $h_x$ has a kernel of dimension 1 and is positive on $T_xM$.  As an  example, by definition characteristic 
hypersurfaces of the D'Alembertian operator on $L$ are lightlike hypersurfaces. Also, horizons (in particular of black holes if any) of subsets of $L$ are topological hypersurfaces and are lightlike whence they are smooth.  

Now we give two opposite classes of examples of lightlike metrics which correspond, by duality to the integrable and contact cases  of the sub-Riemannian situation, respectively:

\subsubsection*{Transversally Riemmannian lightlike metrics} A lightlike metric on a manifold  $I$ of dimension 1 is just 0. Consider now a direct product of $(I, 0)$ with a Riemannian metric $(Q, h)$. This gives a lightlike metric $ h \oplus 0$ on $Q \times I$. A lightlike metric $g$ on  a manifold $M$ is called {\bf transversally Riemannian}Ê  if it is locally isometric to such a (direct) product. This is equivalent to the fact that the flow of any vector field tangent to the characteristic direction $N$ preserves $g$ (it suffices that this happens for one non-singular such vector field). 

\subsubsection*{Generic lightlike metrics} Let $X$ be a vector field tangent to $N$. One sees
that $N$ annualizes the Lie derivative $L_Xg$ (i.e.  $L_Xg(u, v) = 0$, if $u \in N$). Furthermore, $L_Xg$ is conformally well defined: if $X^\prime$ is another vector field tangent to $N$, then $L_{X^\prime}g = f L_Xg$, for some function $f$ on $M$. We say that $g$ is {\bf generic}  if $L_Xg$ has maximal rank, i.e. 
its kernel is exactly $N$. This therefore defines a conformal pseudo-Riemannian structure 
on $TM/N$.

%\subsection{Examples}

\subsection*{Rigidity flavours}  Exactly as in the sub-Riemannian case, lightlike metrics have infinite type and thus are not rigid. Indeed, the local isometry group has infinite dimension 
for  any transversally-Riemannian metric. For example, if $M = \R^{n-1} \times \R$, with 
the  metric $dx_1^2 + \ldots dx_{n-1}^2$, then any map $f(x, t) =  (x,  l(x, t))$ is isometric. 

The key observation of \cite{BFZ} was a kind of Liouville theorem for the lightcone $Co^n$, $n \geq 3$. This is $\R^+ \times S^{n-1}$ endowed with the lightlike  metric $g_{(t, x)}= e^{2t} Can_x$, where $Can$ is the usual metric on $S^{n-1}$. This is in fact the lightcone at 0 in the Minkowski space $Min^{n+1}$. The lorentz group $O^+(1, n)$ acts isometrically 
on $Min^{n+1}$ and hence on $Co^n$. The observation is that any local isometry of $Co^n$ 
coincides with the restriction of an element of $O^+(1, n)$.

%\subsection*{Rigidity of $G$-structures}

%\subsection{Non-rigidity}

%\subsubsection{Foliawithtions}

%\subsubsection{Transversally Riemannian flows}

\subsubsection*{A correspondence:  ``generic transversally conformal Lightlike geometry'' $\leftrightarrow$ Conformal Riemannian geometry } The cone situation generalizes to that 
of {\bf transversally conformal} lightlike structure. This means that the flow of any $X$ tangent to $N$ is conformal, equivalently $L_Xg = fg$ for some function $f$. Locally, $M = Q \times I$ where $I$ is an interval, and $g_{(q, r)} = c(q, r)h_q$, where $h$ is a Riemannian metric on $Q$. 

Assume $\phi$ is an isometry of $(M, g)$, then it acts on $Q$,  the quotient space of its characteristic foliation and induces a diffeomorphism $\psi$, which is obviously conformal for $(Q, h)$. 

Conversely, let $\psi$ a conformal transformation of $(N, h)$, and let us look for an isometry of $(M, g)$ of the form $\phi: (q, r) \to (\phi(q), \delta(q, r))$. We assume here $g$ generic, which means that $\frac{\partial c(q, r)}{\partial r } \neq 0$. Let us assume that $I = \R$, and for any $q$, the map $r \to c(q, r)$ is a global diffeomorphism of $\R$. If $f$ is the conformal distortion of $\psi$, that is $\psi^*h = fh$, then $\phi$
is isometric iff $c(\phi(q), \delta(q, r)) f(q) = c(q, r)$. Our hypotheses imply that for any fixed $q$, $\delta(q, r)$ can be uniquely chosen, and hence a conformal transformation 
of $(N, h) $ admits a unique isometric lifting on $(M, g)$.
(One may compare with a somehow similar construction in  \cite{Fef}).

% generalizes as follows..

%Local expression...

%{\bf generic lightlike}  metric

\section{Sub-rigidity of geometric structures, Results}

%Our principal  (positive) result is the following (where $jet_x^k(\phi) = 1$, means that 
%$\phi$ has the same $k$-jet as the identity at $x$):

We have  the following infinitesimal result for lightlike metrics (where $jet_x^k(\phi) = 1$, means that 
$\phi$ has the same $k$-jet as the identity at $x$):

\begin{theorem} \label{subrigidity} Let $g$ be a generic lightlike metric on a manifold $M$ of dimension $n \geq 4$. Then, a 3-infinitesimal isometry with a trivial 1-jet, has a trivial 2-jet: $$\phi \in \Is^3_x, \; jet^1_x(\phi) = 1 \Longrightarrow jet^2_x(\phi) = 1$$

%A similar statement is true for contact sub-Riemannian metrics on manifolds of dimension 
%$\geq 3$. 

\end{theorem}

This notion was actually brought out by Benveniste-Fisher in \cite{BF} under the name of almost-rigidity. We believe here that the word ``sub-rigid'' is more telling (see also \cite{Dum}). 

In order to keep an elementary level of exposition, we restrict ourselves to geometric structures of order 1,  that is,  $\GL_n(\R)$-equivariant maps $P(M) (= \GL^{(1)}(M)) \to Z$, where $P(M)$ is the frame bundle of $M$ and $Z$ is 
a manifold with a  $\GL_n(\R)$-action.  The classical case of an $H$ structure corresponds to the homogeneous  space $Z = \GL_n(\R)/H$.

\begin{definition} A geometric structure $\sigma$ is $(k+s, k)$-sub-rigid, if any $(k+s)$-isometry 
whose $k$-jet is trivial has a trivial $(k+1)$-jet; formally,  if $\sf{Im}^{k+s, k+1}_x$ denotes the image of $\Is_x^{k+s} \to \Is^{k+1}_x$, then, for any $x$, $\sf{Im}^{k+s, k+1}_x \to \Is^k_x$ is injective.

%that is the  restriction 

\end{definition}

\begin{remarks} ${}$ {\em

1.  $(k+1, k)$-sub-rigidity  means $k$-rigidity.

2.   In particular, a $(k+s, k)$-sub-rigid   structure is  Iso-rigid at order $k$.

3. The theorem above states that generic lightlike metrics 
%and contact sub-Riemannian metrics 
are 
$(3, 1)$-sub-rigid.

}
\end{remarks}

%\bigskip

In the sub-Riemannian case, we have

\begin{theorem}\label{subrigidity.subRiemann} A contact sub-Riemannian metric is 
$(4, 1)$-sub-rigid.

\end{theorem}

\subsubsection*{Example} The paradigmatic example of sub-rigid structures   presented  in \cite{BF} was that of a degenerate framing. That is, on 
$\R^n$, a system of vector fields $x \to (X_1(x), \ldots, X_n(x))$, which are linearly independent everywhere except  at an isolated point, say $0$. As an example, take $n = 1$, and
the geometric structure being a vector field $X(x) = f(x) \frac{\partial}{\partial x }$. A diffeomorphism $\phi$ is isometric if $ \phi^{\prime}(x) f(x) = f(\phi(x))$. If $f$ does not vanish, then we have a true parallelism, and it is $0$-rigid: trivial $0$-jet implies trivial $1$-jet, say at the point $ 0 \in \R$;  in other words, $\phi(0) = 0$ implies $\phi^\prime(0) = 1$. 

Assume now that  $f(0) = 0$, then $\phi$ is isometric up to order $k+1$ at 0, if it satisfies, at the point 0,  all the equations obtained by taking derivatives up to order $k$ of the equality: $ \phi^{\prime}(x) f(x) = f(\phi(x))$. 
Assume $f$ has a zero of order $d$ at 0, e.g. $f(x) = x^d$, and that $jet_0^1(\phi) = 1$,  i.e. $\phi(0) = 0$ and $\phi^\prime(0) = 1$;  then we need derivatives of $\phi$
up to order $d+2$ in order to conclude that $jet^2_0(\phi) = 1$, i.e. $\phi^{\prime \prime}(0) = 0$.  Therefore, the structure is $(d+2, 1)$-sub-rigid.

\begin{remarks} {\em 
${}$

1. An essentially equivalent    example is given  in (\cite{D-G}, \S 5.11.B) to show weakness of    Iso-rigidity  in comparison with rigidity. 

2. One may think following \cite{BF} that, as above,  there is always a degeneracy phenomenon behind sub-rigidity. One may in particular ask if a sub-rigid structure is rigid on an open dense set? However, the examples of the lightlike structure on the Minkowski lightcone, and the standard contact sub-Riemannian metric on the Heisenberg group, show that sub-rigid structures can be homogeneous. They are in particular nowhere rigid. }Ê

\end{remarks}

\section{Proof of Theorem \ref{subrigidity.subRiemann}}

 \begin{proof} Let $(D, h)$ be a contact sub-Riemannian structure on $M$. For computation, it is useful to see $(D, h)$ as an equivalence class of pairs $(\omega, h)$ where $\omega $ is any contact form defining $D$. The correspondence $(\omega, h) \to \bar{h}$ discussed in \S \ref{h.bar} does not depend on the particular choice of $\omega$. We will show that the $1$-jet of $\bar{h}$ is determined by the the $3$-jet of $(\omega, h)$. 
 
 %This would imply that 
 %if $f$ is a $4$-isometry of $(\omega, h)$ is a $2$-isometry for $\bar{h}$, and hence, that 
 %$(h, \omega)$ is $(4, 1)$-rigid (since Riemannian metrics are $(2, 1)$-rigid).
 
 Let $X_1, \ldots X_{2d}$ a local  system of smooth vector fields generating $D$, where $\dim M = 2d +1$. 
 The normalized form
 $\omega^\prime = f \omega$ is defined by $$f d\omega^d(X_1, \ldots, X_{2n}) = \det (h(X_i, X_j)_{ij})$$ 
 
Its  Reeb vector field $R$   is 
  defined algebraically by $$df (X_i) \omega (R) + f d\omega(X_i, R) = 0, \;  \hbox{and}  \; f\omega (R) = 1$$

  If $X_0$ is a vector field transverse to $D$, then  $\bar{h}(X_i, X_0)_{i j}$, $ i, j \geq 0$,  are 
  given by the $h(X_i, X_j)$, $i, j >0$, and the coordinates of $R$ in the moving frame $\{X_0, \ldots X_{2d}\}$. In particular the first derivatives of $\bar{h}$ come from third derivatives of $(\omega, h)$.

 Let $f$ be a diffeomorphism,  $f^*\omega = \omega_1$
 and  $f^*h = h_1$. Thus $f^*\bar{h} = \bar{h_1}$. 
If $f$ is  a $4$-isometry for $(\omega, h)$, then by definition $\omega$ and $\omega_1$ (resp. 
$h$ and $h_1$) coincide up to order 3. It then follows that $\bar{h}$ and $\bar{h_1}$ coincide up to order 1, that is $f$ a 2-isometry for $\bar{h}$. If $f$ was merely a 4-isometry for $(D, h)$, then, $\omega$ will coincide up to order 3 with a multiple $g\omega_1$, which leads to the same conclusion for $\bar{h}$.

   To prove $(4, 2)$-sub-rigidity for $(D, h)$, apply the 1-rigidity (say the $(2, 1)$-sub-rigidity) of Riemannian metrics. We get here that, if  $f$ is a $4$-isometry with a trivial 1-jet, then it has a trivial 2-jet.
   
   % and  has a trivial 2-jet.

 \end{proof}
 
 \subsubsection*{Example } Endow $\R^3$ with the contact  hyperplane field determined by the form  $\omega=dz-xdy$ together  the restriction of $dx^2+dy^2$ on it.  The map $f(x,y,z)=(x+\frac{1}{2}z^2,y-\frac{1}{2}zx^2,z+\frac{1}{2}yz^2)$ belongs to $\Is^3_0$. It  has a trivial 1-jet, but a  non-trivial 2-jet. Therefore, the structure is not $(3, 1)$-sub-rigid.

\section{On the proof of Theorem \ref{subrigidity}}
We give in what follows hints on the proof of Theorem \ref{subrigidity}.
Details,  especially for  \S \S \ref{second.step} and \ref{generalized}
  will appear in \cite{BZ}.

% in the lightlike case.

\subsection{The transversally conformal  case} Let us consider first the transversally conformal case. Locally, $M = Q \times I$, and $g_{(q, r)} = c(q, r)h_q$.  An isometry $\phi$ has the form $\phi: (q, r) \to (\psi(q, r), \delta(q, r))$.  Since $\phi$
preserves the characteristic foliation of $g$, it acts on the quotient space $Q$, that is, $\psi$ does not depend on $r$. 

\subsubsection{The local isometry equation} The isometric equation is:  
$$c(\psi(q), \delta(q, r)) \phi^*h = c(q, r) h$$
That is, $\phi $ is conformal, say with a distortion $f$ (i.e. $\phi^*h = fh$)  which satisfies the  cocycle property: $c(\psi(q), \delta(q, r)) f(q) = c(q, r)$.

%$c(\phi(q), \delta(q, r)) \phi^*h = c(q, r) h$

\subsubsection*{\sc Infinitesimal case} For the sake of  simplicity, even for infinitesimal isometries,  we will assume that their $\psi$-part depends only on $q$ (what is a priori true in this case is that 
the derivatives  on $\psi$ with respect to $r$ vanish according to the order of the infinitesimal isometry).
 
 Fix a point, say $(q, r) = (0, 0)$. The fact that  $\phi= (\psi, \delta)$ is isometric of order 1 (and fixes $(0, 0)$)  means exactly that the previous equation is satisfied for
$(0, 0)$. So $\psi(0) = 0$, $\delta (0, 0) = 0$, and $\psi$ is  conformal at order 1.

\subsubsection*{First step: $\phi \in \Is^3_{(0, 0)(g)}$ and $jet^1_{(0,0)} \phi = 1  \Longrightarrow jet^3_0 \psi = 1$:}

${}$

-- The fact that  $\phi$  has a trivial 1-jet translates to: $D_0\phi = 1$, and $\frac{\partial \delta}{\partial q} = 0$, and $\frac{\partial \delta}{\partial r} = 1$. 

-- The fact that $\phi \in \Is^2_{(0,0)}(g)$ means that  the second derivatives of $\phi$ at $(0, 0)$ satisfy all the equalities obtained by derivating the previous equation. Here, using that
$\frac{\partial \delta}{\partial q} = 0$, we  
   observe that we have in fact   that $\psi \in \Is^2_0(h)$, that is $\psi$ is not only conformal, but isometric for $h$. Then, we use 1-rigidity of 
Riemannian metrics to deduce that $jet^2_0(\psi) = 1$. 

-- Now,   $\phi \in \Is^3_{(0, 0)}(g)$ implies    in particular that $\psi$ is 3-conformal for $h$. We then apply Liouville Theorem, that is the 3-rigidity of conformal Riemannian metrics, and   deduce that $jet^3_0(\psi) = 1$.

\subsubsection{Second step: $\phi \in \Is^3_{(0, 0)}$ and $jet^3_0 \psi = 1  \Longrightarrow jet^2_{(0,0)} \delta = jet^2_{(0, 0)}{r}$, i.e. all the second  derivatives of $\delta$ vanish at $(0, 0)$} \label{second.step}

${}$

-- The equation $\phi \in \Is^3_{(0, 0)}$ obtained by taking second derivatives of the isometric equation gives relations between $jet^2_0(\delta)$ and $jet^3_0(\psi)$. Since, we already know that $D_0\psi = 1$ and all other derivatives of order $\leq 3$ vanish, we get 
equations relating second derivatives of $\delta$ (the first derivatives of $\delta$ are known).  We  then prove that this system of algebraic linear equations  (on these derivatives) is determined and that all 
the second derivatives of $\delta$ vanish.

%The $\phi$-part corresponds  to the action o

\subsection{The general case, generalized conformal structures}  \label{generalized} When $g_{q, r}$ has  a general 
form rather the split one in the transversally conformal case, we get on the quotient space a kind of {\bf generalized conformal structure}. 
This means that at each $q \in Q$, we are giving $S_q  \subset {Sym^2}^*(T_qM)$, the space of Euclidean scalar   products on $T_qQ$, such that $S_q$ is the image of a (non-parameterized) curve. The case of  Riemannian metrics corresponds to   $S_q$ reduced to one point, and that of conformal structures to that where all the elements of $S_q$ are proportional. (Of course, we assume everything depends smoothly on $q$).

The proof of Theorem \ref{subrigidity} goes through  an adaptation of Liouville theorem to generalized conformal structures, that is a generalized conformal structure is $2$-rigid. 
 
 The second step in the proof of Theorem  \ref{subrigidity} is the same as in the transversally conformal case.

\section{Weakness}
% of sub-rigidity}

We show in what follows how the sub-rigidity  is weak in comparison to rigidity. 

\subsection{
Gromov   representation theorem for rigid structures}

Let $G$ be a Lie group acting on a compact manifold $M$ by preserving an analytic  geometric 
structure $\sigma$ and  a volume form. The Gromov representation theorem  concerns the case where $\sigma$ is rigid and $G$ is a simple Lie group. It states that $M$ tends to look like a quotient $G/\Lambda$, where $\Lambda $ is a co-compact lattice in $G$. The precise statement is that $\pi_1(M)$ is large, in the sense that it has a representation in some linear group 
whose the Zariski closure of the image contains a copy of the Lie group $G$. This result was generalized  for actions of lattices in $G$,  by Fisher-Zimmer:   

\begin{theorem} \cite{FZ} Let $\Gamma$ be a lattice in a simple Lie group $G$ of rank $\geq 2$. Suppose  $\Gamma$ acts  on a manifold $M$ analytically by preserving an analytic  rigid geometric  
structure,  and  ergodically for a volume form. Then, either:

1. $\Gamma$ acts via a homomorphism in a compact subgroup $K \subset \Diff(M)$ (and thus 
$M$ is a homogeneous space $K/C$, by ergodicity), or

2. As in Gromov representation theorem, $\pi_1(M)$ admits a homomorphism in some $\GL_N(\R)$
whose Zariski closure   contains a subgroup locally isomorphic to the Lie group $G$.

\end{theorem}
 
We will show below that this does not extend to sub-rigid structures. 

\subsubsection*{Extension of Killing fields} One crucial ingredient in the proof 
on the previous results  is that, for  rigid geometric structures,  local analytic Killing fields extend to the full manifold if it is simply connected (see \cite{Amo, D-G, Gro}). More precisely, let $M$ be    analytic, simply connected  and  endowed with an 
analytic rigid geometric structure $\sigma$. Let $V$ be a Killing field of $\sigma$ defined on an open set $U \subset M$ (that is the local flow of $V$ preserves $\sigma$). Then, $V$ 
extends (as an analytic Killing field) to $M$.  

This fact is no longer true for sub-rigid structures. Indeed, let $g = x^{2d}dx^2$ be a
 ``singular'' Riemannian metric on $\R$. On $\R-\{0\}$, the metric is regular, and hence flat, it has a Lie algebra of Killing fields of dimension 1. No such  Killing field extends at 0. Indeed, as 0 is  the unique singular point of $g$, it will be fixed by any local isometry $\phi$ defined on its neighbourhood. One then shows that   $\phi $ is necessarily $\pm Id$. Indeed, there is a well defined distance $d_g$ derived from $g$. Thus, $d_g(0, x) = d_g(0, \phi(x))$.

\subsection{No Gromov   representation  for sub-rigid structures}

In the sub-rigid case, we have the following example:

\begin{theorem} \label{example} The lattice $\SL_3(\Z)$ acts analytically and ergodically on a compact simply connected manifold, by  preserving  an 
analytic sub-rigid structure and a volume form.  More precisely, there exists a holomorphic action  of $\SL_3(\Z + j\Z)$, $j = e^{\frac{2\pi}{3}i}$,  on a compact  Calabi-Yau 3-manifold (i.e. a simply connected  K\"ahler manifold
with a holomorphic volume form). The action preserves a holomorphic sub-rigid structure, and is ergodic (it is measurably isomorphic to an affine action of a complex torus of dimension 3).

 \begin{remark} Observe this is a   sub-rigid  counter-example  for the Fisher-Zimmer version concerning  higher rank lattices actions. The original Gromov's theorem  deals with  actions  of Lie groups.

 \end{remark}

\end{theorem}

%\subsubsection*{Pullback of geometric structures}

\begin{proof}
Before giving the construction, let us discuss somewhat the general question of taking pull-backs 
 % there is no way to define inverse images 
 of geometric structures.

 \subsubsection*{Pull-Back} 
 Let $\pi: M^\prime \to M$ be a differentiable map, with $M$ and $M^\prime$ of same dimension $n$. Assume $M$ is endowed with an $H$-structure $\sigma$ ($H $ a subgroup of $\GL_n(\R)$). If $\pi$ has no critical points ($\pi$ a local diffeomorphism) then, one defines straightforwardly 
 $\pi^*(\sigma)$. Indeed, $jet^1(\pi)$ is well defined as a map $\GL^{(1)}(M^\prime) \to \GL^{(1)}(M)$, and then one composes it with $\sigma: \GL^{(1)}(M) \to \GL_n(\R)/H$. 
 In contrast, there is generally no mean to define $\pi^*(\sigma)$ at  a critical point. As   an example,  there is no definition  of the inverse image of a vector field on $M$.

 Let us now describe a situation   where the definition of the pull-back of an $H$-structure is possible as a geometric structure in the Gromov sense, 
 but  not as an $H$-structure.   Assume $\sigma$ is   a parallelism $x \to (e_1(x), \ldots, e_n(x))$  on $M$. One defines a geometric structure   $  \sigma^\prime: \GL^{(1)}(M^\prime) \to    \sf{Mat}_n(\R)$ by: 
$$u = (u_1, \ldots, u_n) \in \GL^{(1)}_y(M^\prime) \to \sigma^\prime (u) = (a_{ij}(u)) \in  \sf{Mat}_n(\R)$$
where $D_y\pi (u_i) = \Sigma a_{ij}(u)e_j(\pi(y))$. In other words,    $\sigma^\prime(u)$ is the matrix of the derivative $D_y\pi$ with respect to the bases $(u_1, \ldots, u_n)$ 
and $(e_1(\pi(y)), \ldots, e_n(\pi(y)))$
of $T_yM^\prime$
and  $T_{\pi(y)}M$, respectively.

\subsubsection*{Case of the affine flat connection} Another situation which serves in the proof of our theorem is that of the usual affine structure on $\R^n$. This an $H$-structure of order 2, i.e. a map $\sigma: \GL^{(2)}(\R^n) \to \GL^{(2)}_0(\R^n)$, where $\GL_0^{(2)}(\R^n)$ is the set of invertible jets in $jet^2_0(\R^n,  \R^n)$ (the space of jets $(\R^n, 0) \to \R^n$),  i.e. the inverse image of $\GL_n(\R)$ under the projection $jet^2_{0,0}(\R^n, \R^n) \to 
\GL_n(\R^n) \subset jet^1_{0, 0}(\R^n,  \R^n)$ (the space of jets $(\R^n, 0) \to (\R^n, 0)$). The usual affine connection on $\R^n$ 
is obtained from the projection $\GL^{(2)}(\R^n) = \R^n \times  \GL^{(2)}_0(\R^n) \to \GL^{(2)}_0(\R^n)$. Now the point is that $\sigma$ extends as a map $\bar{\sigma}: jet^2_0(\R^n,  \R^n) \to 
jet^2_{0, 0}((\R^n, \R^n) $. The smooth map  $\pi $ always induces a map
$jet^2(\pi): \GL^{(2)}(M^\prime) \to jet^2_0(\R^n, \R^n)$. We define $\pi^*(\sigma)$
as $\bar{\sigma} \circ jet^2(\pi)$. 

%(regarding notation, the subscribe $0$ in $jet$ means jets of maps sending 0 to any point in %the target, and $0, 0$ means  that the map sends 0 to 0, in $\R^n$).

If the degeneracy of $\pi$ is bounded, that is there exists $k$ such that 
$jet^k_x(\pi) \neq 0$, for any $x$, then $\pi^*(\sigma)$ is sub-rigid.

The construction is natural, and thus, if a group $\Gamma$ acts on $M^\prime$ and $\R^n$
equivariantly with respect to $\pi$, and if the action on $\R^n$ is affine, then, $\Gamma$ preserves the pull-back of the geometric structure on $M^\prime$.

All this applies identically to the torus $\T^n = \R^n/\Z^n$, since we have the same trivialization of the jet bundle $\GL^{(2)}(\T^n)$. 

\subsubsection*{Actions} Consider $\Gamma$ a   subgroup of $\SL_n(\Z)$ and let it act as usually on    $\T^n$. Blow-up a finite orbit of $\Gamma$ (e.g. a rational point), and get a manifold $M^\prime$ with a $\Gamma$-action. It was proved by Katok-Lewis \cite{KL} that this action is volume preserving, and by Beneviste-Fisher \cite{BF} that it preserves a sub-rigid structure, but can not preserve a rigid one.

\subsubsection*{Orbifolds} Here, we assume that there is a finite index subgroup $F \subset \SL_n(\Z)$ commuting with $\Gamma$. We then consider the orbifold $M_0= \T^n/F$. It inherits  a $\Gamma$-invariant natural flat affine connection in an orbifold sense. 

The next step is to desingularize $M_0$ in order to get  a (regular)  manifold $M^\prime$ with 
a $\Gamma$-action equivariant with respect to a projection $M^\prime \to M_0$.

\subsubsection*{An example}
% in dimension 2} 
Take $F$ to be the group   isomorphic 
to $Z/2 \Z$ generated by the involution $I: x \to -x$. If $n = 2$, the quotient around fixed points is just a cone with opening angle $\pi$. It follows that the  so obtained orbifold
is a topological surface.  It is in fact a topological sphere, with exactly 4 conic singularities on which $\Gamma = \SL_2(\Z)$ acts continuously by preserving a continuous volume form. The singularities can be solved to give the usual differentiable structure on 
the sphere, but this can not be done $\Gamma$-equivariantly. 

%This is not anymore true 

 %In dimension $n >2$, $M_0$ is    anymore a topological manifold. One blow-up $\T^n$ on the %fixed points of $I$ and obtain a space $M_1$. The involution $I$ acts on $M_1$, and there,    %singularity are conical, and hence $M_2 = M_1/F$ is a topological manifold. In this case, one %can revolve singularities in a $\Gamma$-equivariant meaner.   

\subsubsection*{Complex case}  A higher dimensional generalization is possible, but in a complex framework. So, we start with a  complex torus  $A = \C^n/ \Lambda$ (of complex dimension $n$). We consider $M_0 = A/F$, where $F$ is the same previous group generated by the involution $I$. If $n \geq 2$, $M_0$ is no longer a topological manifold, since the fixed points of $I$ are not conical. We then start blowing-up $A$ on $I$-fixed points, in a complex way, and get $M_1$. We have a projection $p_1: M_1 \to A$ with singular fibers  isomorphic to $\C P^{n-1}$ over $I$-fixed points. Now, $F$ acts naturally on $M_1$ with conical singularities, and hence $M_2 = M_1/F$ is a topological manifold. The resolution of singularities yields a complex structure on $M_2$ with a  natural $\Gamma$-holomorphic action. 

\subsubsection*{Our case} For $n=2$, we get a Kummer surface, a special case of K3 surfaces. Observe that the volume form $ dv= dz_1 \wedge dz_2$ is $F$-invariant and hence well defined on $M_2$. However, even if the form $p_1^*(dv)$ is  singular along the exceptional fibers, it gives rise to a true regular holomorphic volume form on $M_2$.

In order to have a similar construction in dimension 3, we replace $F$ by the   group generated by the rotation $J: z \to jz$  where $j = e^{\frac{2 \pi}{3}i}$ on $\C^3$ (see \cite{CZ}, \S 7.6). It preserves the volume form, and therefore, we get on the corresponding $M_2$ a holomorphic volume form.

Regarding the   $\Gamma$-action, we take, $\Gamma = \SL_3(\Z + j\Z)$ and $\Lambda = (\Z + j\Z)^3 \subset \C^3$. 
  Thus,  $\Gamma$   is a lattice in $\SL_3(\C)$, it preserves $\Lambda$  and commutes with $J$.  
  
  As in the cases  $n =1, 2$, one can prove directly that $M_2$ is simply connected.  Another idea is to use the fact that $M_2$ has  holomorphic volume form to deduce it has a vanishing first Chern class. Then, apply Yau's theorem to get a K\" ahler Ricci flat metric on it. But, for such manifolds, up to a finite cover, there is  a de Rham decomposition into a product of a flat torus, and (compact) simply connected manifolds (hyper-K\"ahler and Calabi-Yau, see \cite{Bea}). Thus, it suffices to verify  that $M_2$ has no torus (of dimension 1, 2 or 3) as a factor, to prove that it has a finite fundamental group. 
\end{proof}
 
\begin{remark} {\em One can use general theory of rigid transformation groups to see that 
the latter  action can not preserve a (real) analytic rigid geometric structure. Indeed, by \cite{D-G, Gro}, the isometry group of a unimodular  analytic rigid structure on a simply connected manifold have a finite number of connected components. This means that up to a finite index, the $\Gamma$-action extends to an action of a Lie group, which  can be easily seen to be impossible. } 

\end{remark}

\end{document}